\documentclass[11pt]{article}
\usepackage[utf8]{inputenc}
\usepackage{geometry}
\usepackage{amsmath}
\usepackage{amsfonts}
\usepackage{amssymb}
\usepackage{amsthm}
\usepackage{graphicx}
\usepackage{enumitem}
\usepackage{bm}
\usepackage{caption}
\usepackage{subcaption}
\usepackage{url}
\usepackage[table,xcdraw]{xcolor}
\usepackage{mathtools}
\usepackage{systeme}
\usepackage{hyperref}
\usepackage[capitalise]{cleveref}
\usepackage{thm-restate}

\newtheorem{theorem}{Theorem}[section]

\newtheorem{proposition}[theorem]{Proposition}

\newtheorem{lemma}[theorem]{Lemma}

\newtheorem*{remarks*}{Remarks}
\newtheorem*{remark*}{Remark}

\crefname{equation}{}{}
\crefname{item}{}{}
\numberwithin{equation}{section}

\newcommand{\N}{\mathbb{N}}

\newcommand{\X}{\mathcal{X}}

\def\ew#1{}
	\renewcommand{\ew}[1]{\footnote{\textbf{EW:}#1}}

\title{A proof of a conjecture of Erd\H{o}s and Gy\'{a}rf\'{a}s on monochromatic path covers}

\author{Alexey Pokrovskiy\thanks{Department of Mathematics, University College London, Gower Street, WC1E 6BT, UK. Emails: \texttt{a.pokrovskiy@ucl.ac.uk, lversteegen.math@gmail.com, ella.williams.23@ucl.ac.uk.}} \and Leo Versteegen\footnotemark[1] \and Ella Williams\footnotemark[1]}

\date{}

\begin{document}

\setlength{\baselineskip}{15pt}
\maketitle

\begin{abstract}
In 1995, Erd\H{o}s and Gy\'{a}rf\'{a}s proved that in every $2$-edge-coloured complete graph on $n$ vertices, there exists a collection of $2\sqrt{n}$ monochromatic paths, all of the same colour, which cover the entire vertex set. They conjectured that it is possible to replace $2\sqrt{n}$ by $\sqrt{n}$. We prove this to be true for all sufficiently large~$n$.
\end{abstract}

\section{Introduction}
In \cite{gerencser1967ramsey}, Gerencs{\'e}r and Gy\'{a}rf\'{a}s observed the following fact with a very short and elegant proof.

\begin{theorem}\label{thm:gyarfs-og}
    The vertex set of every 2-edge-coloured complete graph on $n$ vertices can be covered by two monochromatic paths.
\end{theorem}

Importantly, the two monochromatic paths do \emph{not} have to be of the same colour. In his survey \cite{GYARFAS2016}, Gy\'{a}rf\'{a}s mentions that when he first told Erd\H{o}s about this result, Erd\H{o}s misunderstood this detail, which led the two of them to investigate the smallest number of monochromatic paths of the same colour that cover all vertices of a 2-edge-coloured complete graph. In \cite{erdos1995vertex}, they proved the following.

\begin{theorem}\label{thm:erdos-gyarfas}
The vertex set of every 2-edge-coloured complete graph on $n$ vertices can be covered by $2\sqrt{n}$ monochromatic paths, all of the same colour.
\end{theorem}

We remark that although \cref{thm:gyarfs-og} still holds if we additionally impose that the two paths are vertex-disjoint, this is not the case for \cref{thm:erdos-gyarfas}, where paths are allowed to intersect.

The bound $2\sqrt{n}$ in \Cref{thm:erdos-gyarfas} is best possible up to a constant by the following construction. Suppose $\sqrt{n} \in \mathbb{N}$ and partition $V(K_n)$ into disjoint sets $A$ and $B$ of order $n-\sqrt{n}+1$ and $\sqrt{n}-1$ respectively. Colour the edges contained in $A$ blue, and all other edges red. Each red path covers at most $\sqrt{n}$ vertices in $A$ since it must alternate between $A$ and $B$, so we need at least $\lceil |A|/\sqrt{n}\rceil = \sqrt{n}$ red paths to cover the entirety of $V(K_n)$. Every collection of blue paths covering $V(K_n)$ must include all individual vertices in $B$ as distinct paths and needs another blue path to cover $A$, thus has size at least $|B|+1 = \sqrt{n}$. If $n$ is not a perfect square, one can use the same construction with $\vert B\vert=\lfloor \sqrt{n} \rfloor-1$, so that any cover requires at least $\lfloor \sqrt{n} \rfloor$ paths of the same colour.

Erd\H{o}s and Gy\'{a}rf\'{a}s conjectured (see \cite{GYARFAS2016}) that this construction is the colouring that requires the most paths for a cover, meaning that the bound in their theorem could be improved to $\sqrt{n}$. In this paper, we confirm their conjecture for all sufficiently large $n$.

\begin{theorem}\label{thm:main}
For all $n> 20^{40}$, the vertex set of every 2-edge-coloured complete graph on $n$ vertices can be covered by $\sqrt{n}$ monochromatic paths, all of the same colour.
\end{theorem}

We do not attempt to optimise the constant $20^{40}$, but with some additional technical effort, one can show that for all $n\in \N$, $\sqrt{n}+10$ monochromatic paths of the same colour are sufficient to cover $V(K_n)$.

Several related problems of this type have also been studied, looking for monochromatic structures that either cover or partition the vertex set within an edge-coloured graph.  
For example, Lehel's conjecture, proven for large $n$ by Łuczak, Rödl
and Szemer\'edi \cite{Luczak1998PartitioningTC} and independently by Allen \cite{ALLEN_2008}, and for all $n$ by Bessy and Thomass\'e \cite{BESSY2010176}, strengthens \Cref{thm:gyarfs-og}. It states that the vertex set of every $2$-edge-coloured $K_n$ can be covered by two vertex-disjoint monochromatic cycles (of course, they must be allowed to have different colours). 
A generalisation of \cref{thm:gyarfs-og} conjectured by Gy\'{a}rf\'{a}s \cite{gyarfas-mono-paths}, stating that the vertex set of every $r$-edge-coloured complete graph can be covered by $r$ vertex-disjoint monochromatic paths, remains open. We refer the reader to the survey \cite{GYARFAS2016} for partial results and further related problems. 

Finally, as an interpolation between \Cref{thm:gyarfs-og} and \Cref{thm:erdos-gyarfas}, Eugster and Mousset \cite{eugster-mousset} investigated the minimal number of monochromatic paths needed to cover the $r$-edge-coloured $K_n$ if at most $s\leq r$ different colours are allowed to be used for the paths in total.

\section{Preliminaries}\label{sec:preliminaries}
 We call a collection $\mathcal{P}$ of paths in a red-blue edge-coloured graph $G$ a \textit{red path cover} if every vertex in $V(G)$ is covered by some path in $\mathcal{P}$, and every $P \in \mathcal{P}$ contains only red edges. We define a \textit{blue path cover} similarly. 
The length of a path is the number of edges it contains, and importantly we allow paths to have length zero, so that a path can be a single vertex. 
 Note that paths need not be disjoint, unlike those required for some definitions of the path covering number of a graph.

We will use a bipartite Ramsey-type result for paths determined by Gy\'{a}rf\'{a}s and Lehel in \cite{gyarfas1973ramsey} (see Theorem 3 and Remark 1).
\begin{lemma}\label{lem:bipRamsey}
    If $k$ and $\ell$ are distinct natural numbers then every red-blue edge-coloured complete bipartite graph $K_{\left\lceil\frac{k+\ell}{2}\right\rceil,\left\lceil\frac{k+\ell}{2}\right\rceil}$ contains a red path of length $k$ or a blue path of length $\ell$, i.e., a red path on $k+1$ vertices or a blue path on $\ell+1$ vertices. 
\end{lemma}

Suppose we have found a long monochromatic path $P$ in a red-blue edge-coloured $K_n$, without loss of generality we may assume $P$ is blue. 
We may consider the bipartite graph between $V(P)$ and a subset of $V(K_n) \setminus V(P)$, comprised only of red edges. If there is a long red path, say $R$, in this bipartition then we obtain a large set $V(P)\cap V(R)$ covered by both a red path and a blue path, which allows us some flexibility for covering the remainder of the vertices.  
A similar idea may be used if we can find few red paths that cover many vertices in the bipartition. In the following two lemmas we consider certain degree properties within bipartite graphs that allow us to find few paths covering many vertices.

\begin{lemma}\label{lemma:bipartite-path}
Let $G$ be a bipartite graph on vertex set $X\cup Y$ with the property that every vertex in~$Y$ has degree at least $(\vert X\vert +\vert Y\vert)/2$. Then $G$ contains a path covering $2\vert Y\vert$ vertices.
\end{lemma}
\begin{proof}
    Let $P=x_1y_1\ldots x_ky_k$ be a path starting in $X$ and ending in $Y$ and suppose that $k<\vert Y\vert$. Fix $y_{k+1}\in Y\setminus \{y_1,\ldots,y_k\}$. We have
    \begin{equation*}
        \vert N(y_k)\cap N(y_{k+1})\vert \geq \frac{\vert X\vert +\vert Y\vert}{2}-\left(\vert X\vert -\frac{\vert X\vert+ \vert Y\vert}{2}\right)=\vert Y\vert>k.
    \end{equation*}
    Therefore, $y_{k}$ and $y_{k+1}$ have a common neighbour $x_{k+1}\in X\setminus \{x_1,\ldots,x_k\}$, and we can extend $P$ via $y_kx_{k+1}y_{k+1}$.
\end{proof}

\begin{lemma}\label{lemma:bipartite-decomposition}
Let $m\in \N$ and let $G$ be a bipartite graph on vertex set $X\cup Y$ with  $\vert X\vert \geq \vert Y\vert+2m$, and such that every vertex in $Y$ has degree at least $\vert X\vert -m$. Then $G$ contains a collection of at most $ \lfloor \vert X\vert /\vert Y\vert \rfloor$ paths whose union covers all but $\vert Y\vert + 2m$ vertices in $X$ and all vertices in $Y$.
\end{lemma}
\begin{proof}
    We prove the statement by induction on $|X|$, always assuming that $\vert X\vert \geq \vert Y \vert + 2m$. Since
    \begin{equation*}
        \vert X\vert -m \geq \frac{\vert X\vert+\vert Y\vert}{2},
    \end{equation*}
    we may apply \Cref{lemma:bipartite-path}, to obtain a path $P$ covering $2\vert Y\vert$ vertices in $G$. In particular, $P$ covers all vertices in $Y$. Let $X'=X\setminus V(P)$, and note that in the induced subgraph $G'=G[X'\cup Y]$, every vertex in $Y$ has degree at least $\vert X'\vert-m$. For the base case of our induction, if $\vert X'\vert \leq \vert Y\vert+2m$, then $\{P\}$ is a collection of paths as desired. Otherwise, we may apply the induction hypothesis to $G[X'\cup Y]$ to obtain a collection of paths $P_1,\ldots,P_k$ with $k\leq \lfloor \vert X'\vert /\vert Y\vert\rfloor=\lfloor \vert X\vert /\vert Y\vert\rfloor-1$ that cover all but $\vert Y\vert +2m$ vertices in $X'$. Adding $P$ to this collection yields the claim.
\end{proof}

To get the exact bound in \Cref{thm:main}, we also require a version of \Cref{lemma:bipartite-decomposition} which lets us cover all vertices in $X$ when many vertices in $Y$ are adjacent to all vertices in $X$.

\begin{lemma}\label{lemma:bipartite-decomposition-Y0}
Let $G$ be a bipartite graph on vertex set $X\cup Y$ with sets $X_0=\{x\in X: d(x)=\vert Y\vert\}$, $X_1=X\setminus X_0$, $Y_0=\{y\in Y:d(y)=\vert X\vert\}\neq \emptyset$, and $Y_1=Y\setminus Y_0$ such that
\begin{enumerate}[label=\upshape(\roman*)]
    \item $\vert X\vert> \vert Y\vert$, and
    \item $X_1=Y_1=\emptyset$ or $\vert X_0\vert/\vert Y_1\vert>2\vert X_1\vert /\vert Y_0\vert$.
\end{enumerate}
Then $G$ contains a collection of at most $ \left\lceil\frac{ \vert X\vert}{\vert Y\vert +1}\right\rceil$ paths that cover all vertices in $G$.
\end{lemma}
\begin{proof}
    We prove the statement again by induction on $|X|$, always assuming that $|X|\geq |Y|+1$. The base case $|X|=1$ is trivial, because there exists no graph $G$ as above with $1=|X|>|Y|\geq |Y_0|>0$. It follows from (i) and (ii) that $\vert X_0\vert \geq \vert Y_1\vert+1$. Indeed, if $X_1 = Y_1 = \emptyset$ then it is clear, and so supposing otherwise that $\vert X_0\vert/\vert Y_1\vert>2\vert X_1\vert /\vert Y_0\vert$ but that $|X_0|\leq |Y_1|$, then from (i) we must have $|X_1|>|Y_0|$, which cannot all hold true at the same time. Therefore assuming $\vert X_0\vert \geq \vert Y_1\vert+1$, $G$ contains a path~$P$ with $2\vert Y_1\vert+1$ vertices, alternating between $X_0$ and $Y_1$, such that $P$ starts at $x_0 \in X_0$ and also ends in $X_0$. In particular $P$ covers all vertices in $Y_1$. 
    
    Furthermore, since every vertex in $Y_0$ has neighbourhood $X$, and $\vert X\vert > \vert Y\vert$, $G$ contains a path~$Q$ with $2\vert Y_0\vert$ vertices, alternating between $Y_0$ and $X\setminus V(P)$, that starts at $y_0\in Y_0$ and covers $\min\{\vert Y_0\vert, \vert X_1\vert\}$ vertices in $X_1$. This is possible since $|X\setminus V(P)|= |X|-|Y_1|-1\geq |Y|-|Y_1| = |Y_0|$. Since $V(P)\subseteq X_0 \cup Y_1$ and $V(Q)\subseteq (X\setminus V(P))\cup Y_0$, we have that $P$ and $Q$ are vertex-disjoint. Thus we can consider a concatenated path $R=Px_0y_0Q$ and let $X'=X\setminus V(R)$. Note that $|X'|= |X| - |Y|-1$. 
    
    Observe that if $|X'|=0$, i.e., if $|X| = |Y|+1$, then the collection consisting only of the path $R$ covers all vertices in $G$, as required by the lemma. Now, let us assume $|X'|>0$ and let $G'=G[X'\cup Y]$. Defining the sets $X_0',X_1',Y_0'$, and $Y_1'$ analogously for $G'$ to how $X_0,X_1,Y_0$, and $Y_1$ were defined for $G$ respectively, we observe that $Y_1'\subseteq Y_1$ and $X_1'= X_1\setminus V(Q)$, while $Y_0\subseteq Y_0'$. In particular, $Y_0'\neq \emptyset$. If we assume $X_1'=\emptyset$ or $Y_1'=\emptyset$, then it follows that $G'$ is complete bipartite and it is trivially possible to cover $X\setminus V(R)$ with $\left\lceil\frac{ \vert X'\vert}{\vert Y\vert +1}\right\rceil =\left\lceil\frac{ \vert X\vert}{\vert Y\vert +1}\right\rceil-1$ paths so that nothing is left to show. Suppose therefore that $X_1',Y_1'\neq \emptyset$. In particular, we must have $\vert X_1\vert >\vert Y_0\vert$ so that $\vert V(Q)\cap X_1\vert =\vert Y_0\vert$ and $R$ covers all of $Y$. Using the former, we see that

    \begin{equation}\label{eqn:fractionG'}
        \frac{\vert X_0'\vert}{\vert Y_1'\vert}\geq\frac{\vert X_0\setminus V(P)\vert}{\vert Y_1\vert}\geq\frac{\vert X_0\vert}{\vert Y_1\vert}-2>2\left(\frac{\vert X_1\vert}{\vert Y_0\vert} - 1 \right)= \frac{2\vert X_1\setminus V(Q)\vert}{\vert Y_0\vert}\geq\frac{2\vert X_1'\vert}{\vert Y_0'\vert}.
    \end{equation}

    If $\vert X'\vert> \vert Y\vert$, we may use the induction hypothesis to obtain a family $\mathcal{P}$ of at most $\left\lceil\frac{ \vert X'\vert}{\vert Y\vert +1}\right\rceil=\left\lceil\frac{ \vert X\vert}{\vert Y\vert +1}\right\rceil-1$ paths which cover all of $X'$ and $Y$, and together with $R$, these paths cover all of $G$. 
    
    If $\vert X'\vert\leq\vert Y\vert$ on the other hand, then we must have $\vert Y_0'\vert >\vert X_1'\vert$. Indeed if this were not the case, by \cref{eqn:fractionG'} we would have $|X_0'|>2|Y_1'|$ implying $|X'| = |X_0'|+|X_1'|> |Y_0'|+2|Y_1'|\geq |Y|$, a contradiction. So, assuming $\vert Y_0'\vert >\vert X_1'\vert$, it is possible to cover all of $X_1'$ with a single path $Q'$ of length $2\vert X_1'\vert$. The remaining graph $G'-V(Q')$ is complete bipartite on bipartition classes of size $|X'|-|X_1'|$ and $|Y|-|X_1'|$, so using $|X'|-|X_1'|\le |Y|-|X_1'|$ we may find a path $P'$ that covers all of $X_0'$ in it. Concatenating $P'$ and $Q'$ suitably, we obtain a second path $R'$ such that $R$ and $R'$ cover all of $G$.
\end{proof}

\section{Proof of \Cref{thm:main}}

Let $n \in \N$. In what follows, without loss of generality we assume the vertex set of $K_n$ is given by $[n]$. Let $\X_n$ be the set of all red-blue colourings of $E(K_n)$. For $\chi \in \X_n$, denote by $f(n,\chi)$ the size of a smallest family of monochromatic paths of the same colour that cover $[n]$. Furthermore, let $f(n)=\max_{\chi\in \X_n} f(n,\chi)$.

In order to prove \cref{thm:main}, we follow an inductive strategy to initially prove a weaker upper bound on $f(n)$ for all $n \in \N$, given in \cref{prop:main} which can then be bootstrapped to prove our main theorem. Thus in each of the lemmas in this section, we will assume an upper bound on $f(m)$ for all $m<n$ that will relate to our inductive hypothesis.

We first introduce a lemma which tells us that under a colouring $\chi \in \X_n$, if there exists a large set $S$ that is covered by both a collection of few red paths, and a collection of few blue paths, then we can obtain a strong upper bound on $f(n,\chi)$.

\begin{lemma}\label{lemma:reduction}
Let $C_1\geq C_2\geq 0$ and let $n\in \N$ be such that $f(m)<\sqrt{m}+C_1$ for all $m<n$. Let further $k\in \N$ and $\chi\in \X_n$ be such that there exist red paths $P_1,\ldots, P_k\subset K_n$, blue paths $Q_1,\ldots,Q_k\subset K_n$, and a set $S\subset [n]$ satisfying $S\subset (V(P_1)\cup \ldots \cup V(P_k))\cap (V(Q_1)\cup \ldots \cup V(Q_k))$. If
\begin{equation*}
    \sqrt{n-\vert S\vert}+C_1+k\leq \sqrt{n}+C_2,
\end{equation*}
then $f(n,\chi)\leq \sqrt{n}+C_2$. In particular, if there is a red path $P$ and a blue path $Q$ such that $S\subset V(P)\cap V(Q)$, and
\begin{equation*}
    \vert S\vert \geq 2(C_1-C_2+1)\sqrt{n},
\end{equation*}
then $f(n,\chi)\leq \sqrt{n}+C_2$.
\end{lemma}

\begin{proof}
Suppose there exist red paths $P_1,\ldots, P_k\subset K_n$, blue paths $Q_1,\ldots,Q_k\subset K_n$, and a set $S\subset [n]$ satisfying  $S\subset (V(P_1)\cup \ldots \cup V(P_k))\cap (V(Q_1)\cup \ldots \cup V(Q_k))$.
Let $H$ be the 2-edge-coloured graph obtained by deleting $S$ from $K_n$ under the colouring $\chi$. Then $H$ has $m=n-|S|$ vertices, and by definition of $f(m)$, there exists a collection of at most $f(m)$ monochromatic paths of the same colour which cover $V(H)$. Adding $\{P_1,\ldots,P_k\}$ or  $\{Q_1,\ldots,Q_k\}$ to this collection according to whether these paths are red or blue respectively yields a collection of at most $f(m)+k $ monochromatic paths of the same colour covering $V(H) \cup S = [n]$ in $K_n$. By our assumptions, we have 
$$f(n,\chi) \leq f(m) + k < \sqrt{m} + C_1 +k = \sqrt{n-\vert S\vert}+C_1+k\leq \sqrt{n}+C_2,$$
as desired. 
For the second part of the statement, observe that if there is a red path $P$ and a blue path $Q$ such that $S\subset V(P)\cap V(Q)$ for some $S \subset [n]$ with $\vert S\vert \geq 2(C_1-C_2+1)\sqrt{n}$, then
\begin{equation*}
    \sqrt{n-|S|} \leq  \sqrt{n-2(C_1 -C_2+1)\sqrt{n}} \leq \sqrt{n} - (C_1 -C_2+1).
\end{equation*}
Rearranging gives 
\begin{equation*}
    \sqrt{n-\vert S\vert}+C_1+1\leq \sqrt{n}+C_2,
\end{equation*}
and so taking $k=1$, we deduce that $f(n,\chi)\leq \sqrt{n}+C_2$.
\end{proof}

We now show that if $\chi\in \X_n$ is such that there exists no small monochromatic path cover, then there must exist some long monochromatic path, of colour $\gamma$ say, such that for every vertex $y$ outside of the path, there are few vertices $x$ on the path such that $\chi(xy)=\gamma$.

\begin{lemma}\label{lemma:long-path}
Let $C_1\geq C_2\geq 0$ and let $n>10^4(C_1-C_2+1)^4$ be such that $f(m)<\sqrt{m}+C_1$ for all $m<n$. If $\chi\in \X_n$ is such that $f(n,\chi)>\sqrt{n}+C_2$, then there is a monochromatic path $P$ such that the size of $Y= [n]\setminus V(P)$ is at most $\sqrt{n}+10(C_1-C_2+1)\sqrt[4]{n}$. Furthermore, if $\gamma$ is the colour of the edges of $P$, then for every $y\in Y$, the size of $\{x\in V(P): \chi(xy)=\gamma\}$ is at most $2(C_1-C_2+1)\sqrt{n}$.
\end{lemma}

\begin{proof}
Let $\Delta = C_1-C_2+1$ and let $\chi\in \X_n$ be such that $f(n,\chi)>\sqrt{n}+C_2$. By \Cref{thm:gyarfs-og}, without loss of generality, there exists a blue path $Q$ covering $\lfloor n/2\rfloor$ vertices. Let $W$ be a subset of $[n]\setminus V(Q)$ of size $\lfloor n/2\rfloor$ and let $H$ be the complete bipartite graph between $V(Q)$ and $W$ under the colouring $\chi$. Let $R$ be a longest red path in $H$. The set $V(Q)\cap V(R)$ is covered both by a single red and a single blue path so that by \Cref{lemma:reduction}, we may assume $\vert V(Q)\cap V(R)\vert<2\Delta\sqrt{n}$. However, since $R$ is a path in $H$, it must alternate between $V(Q)$ and $W$ and so at least $\lfloor \vert V(R)\vert /2\rfloor$ vertices of $R$ lie in $V(Q)$. 
We can therefore obtain the following
\begin{equation*}
    \frac{|V(R)|-1}{2}\leq \left\lfloor \frac{|V(R)|}{2} \right\rfloor \leq |V(Q)\cap V(R)| < 2\Delta\sqrt{n}\leq \lceil 2\Delta\sqrt{n} \rceil. 
\end{equation*}
Rearranging, we deduce that $|V(R)|<2\lceil2\Delta\sqrt{n}\rceil+1$. Note that $H$ is a red-blue edge-coloured $K_{\lceil (n-1)/2 \rceil, \lceil (n-1)/2 \rceil}$, so applying \Cref{lem:bipRamsey} with $k = 2\lceil2\Delta\sqrt{n}\rceil$ and $\ell = n - 1 -k$ tells us that $H$ contains a red path on $k+1$ vertices or a blue path on $\ell +1$ vertices. Since $R$ was chosen to be a longest red path in $H$, the former case cannot occur. Thus there exists a blue path covering at least $\ell +1 =n-k=n-2\lceil2\Delta\sqrt{n}\rceil \geq n - 5\Delta\sqrt{n}$ vertices in $H$.

Now let $P$ be a longest blue path in $K_n$, denote its vertex set by $X$, and let $Y=[n]\setminus X$. Then by the above argument, $\vert X\vert \geq n-5\Delta\sqrt{n}$ and hence $\vert Y\vert\leq 5\Delta\sqrt{n}$. Consider $y\in Y$ and let $B$ denote the set of blue neighbours of $y$ in $P$. Since $P$ cannot be extended, $B$ does not contain any endpoint or any two consecutive vertices of $P$. Furthermore, if we direct $P$ arbitrarily, the predecessors of $B$ form a red clique set. Indeed, if $P=v_1v_2\ldots v_r$, and $v_iv_j$, $v_{i+1}y$, and $v_{j+1}y$ were all blue for some $1\leq i<j\leq r$, then $v_1\ldots v_iv_jv_{j-1}\ldots v_{i+1}yv_{j+1}\ldots v_r$ would be a blue path covering $r+1$ vertices, contradicting the maximality of $P$. However, we may argue as before that $V(P)$ cannot contain a red path on $2\Delta\sqrt{n}$ vertices, which means that $y$ has at most $2\Delta\sqrt{n}$ blue neighbours in $X$, as required by the ``furthermore" part of the lemma.

In other words, taking $m=\lfloor2\Delta\sqrt{n}\rfloor$, each vertex in $Y$ has at least $\vert X\vert -m$ red neighbours in $X$. Thus, as $n>196\Delta^2$,  $\vert X\vert \geq n-5\Delta\sqrt{n}$, and $\vert Y\vert\leq 5\Delta\sqrt{n}$,
\begin{equation*}
    \vert X\vert \geq n-5\Delta\sqrt{n} > 14\Delta\sqrt{n}-5\Delta\sqrt{n}=9\Delta\sqrt{n}\geq \vert Y \vert + 2m.
\end{equation*}
It now follows from \Cref{lemma:bipartite-decomposition} that the bipartite graph given by the red edges between $X$ and $Y$ contains a collection of at most $\lfloor\vert X\vert/\vert Y\vert\rfloor <n/\vert Y\vert$ paths which cover all but $|Y|+ 2m\leq 9\Delta\sqrt{n}$ vertices in $X$ and all vertices in $Y$.

However, all vertices in $X$ that are covered by this collection of paths are also covered by the single blue path $P$, which means that by \Cref{lemma:reduction},

\begin{equation*}
    \sqrt{n}+C_2<\sqrt{14\Delta\sqrt{n}}+C_1 +n/\vert Y\vert\leq 4\Delta\sqrt[4]{n} +C_1+n/\vert Y\vert.
\end{equation*}
It follows that
   $ \sqrt{n}<n/\vert Y\vert + 5\Delta\sqrt[4]{n},$
which implies
\begin{equation*}
    \vert Y \vert < \frac{n}{\sqrt{n}-5\Delta\sqrt[4]{n}}=\frac{\sqrt{n}}{1-5\Delta/\sqrt[4]{n}}.
\end{equation*}
Since we have $\sqrt[4]{n}\geq 10\Delta$ by assumption and $1/(1-x)\leq1+2x$ for $x\in [0,1/2]$, we obtain
\begin{equation*}
    \vert Y \vert < \sqrt{n}+10\Delta\sqrt[4]{n},
\end{equation*}
as desired.
\end{proof}

If there exists a monochromatic path in $K_n$ containing almost all vertices, we can afford to cover the remaining vertices individually or in pairs. The following lemma formalises this. 

\begin{lemma}\label{lemma:Y-bound}
    Let $n\in \N$ and $\chi\in \X_n$ and suppose $P$ is a blue path in $K_n$. Let $Y=[n]\setminus V(P)$ and write $Y_0$ for the set of vertices in $Y$ that are not contained in a blue edge with any vertex in $V(P)$.
    Then we have
    \begin{equation*}
         f(n,\chi) \leq 1+\left\lceil \frac{\vert Y\setminus Y_0\vert}{2} \right\rceil + \vert Y_0\vert < 2+ \vert Y\vert/2 + \vert Y_0\vert/2.
    \end{equation*}
    
\end{lemma}

\begin{proof}
    Let $Y_1\coloneqq Y\setminus Y_0$ and note that for any two vertices $y_1,y_2\in Y_1$, there exists a blue path between $y_1$ and $y_2$ whose interior points lie in $P$. Thus, by partitioning the vertices in $Y_1$ arbitrarily into pairs, we may cover all of them with at most $\lceil \vert Y_1\vert/2 \rceil$ blue paths. Adding a single-vertex path for each $y\in Y_0$ and $P$ to this collection, we obtain a blue path cover of $K_n$ of size at most $1+\lceil\vert Y_1\vert/2\rceil+\vert Y_0\vert$.
    Since any monochromatic path cover of $[n]$ must have size at least $f(n,\chi)$, this yields
\begin{equation*}
   f(n,\chi) \leq 1+\left\lceil \frac{ \vert Y_1\vert}{2}\right\rceil+ \vert Y_0\vert <  2 + \frac{|Y_1|}{2} + |Y_0| =2 + \frac{|Y|-|Y_0|}{2} + |Y_0| = 2+ \frac{ \vert Y\vert}{2} + \frac{\vert Y_0\vert}{2}. 
\end{equation*}
\end{proof}

We now prove a weak version of \Cref{thm:main}, which will bootstrap to prove \Cref{thm:main} in full.

\begin{proposition}\label{prop:main}
For all $n\in \N$ and $C=20^4$, $f(n)<\sqrt{n}+C$.
\end{proposition}

\begin{proof}[Proof of \Cref{prop:main}] 
We prove the claim by induction on $n$, noting that it is trivial for $n\leq C$. Suppose therefore that $n>C$ and that $f(m)<\sqrt{m}+C$ for all $m<n$ and consider a colouring $\chi\in \X_n$. Applying \Cref{lemma:long-path} with $C=C_1=C_2$, we obtain without loss of generality a blue path $P$ such that the size of $Y=[n]\setminus V(P)$ is at most $\sqrt{n}+10\sqrt[4]{n}\leq 3\sqrt{n}/2$ and such that every $y\in Y$ has at least $\vert X\vert - 2\sqrt{n}$ red neighbours in $X=V(P)$. Let $Y_0$ be the set of vertices in $Y$ that are not contained in a blue edge with any vertex in $V(P)$, so that by \Cref{lemma:Y-bound}, if $\vert Y_0\vert \leq \sqrt{n}/2$ or $\vert Y\vert\leq \sqrt{n}$, then $f(n,\chi)\leq \sqrt{n}+2$ and there is nothing left to show. Let us therefore assume that $\vert Y_0\vert >\sqrt{n}/2$ and $\vert Y\vert > \sqrt{n}$. 

On the other hand, we have that $\vert X \vert \geq n-3\sqrt{n}/2\geq 3\sqrt{n}/2+2m$ for $m=2\sqrt{n}$. Therefore, by applying \Cref{lemma:bipartite-decomposition}, we see that there must exist a collection of at most $\sqrt{n}$ red paths that cover all vertices in $Y$ and all but at most $6\sqrt{n}$ vertices in $X$. Let $X'\subset X$ denote the vertices that are not covered by these paths and recall that all edges between $Y_0$ and $X$, in particular those between $Y_0$ and $X'$, are red. Since $\vert Y_0\vert>\sqrt{n}/2$ and $\vert X'\vert <6\sqrt{n}$ we can cover all of $X'$ with at most 12 red paths between $X'$ and $Y_0$, which completes the proof.
\end{proof}

We are now ready to prove \Cref{thm:main}.

\begin{proof}[Proof of \Cref{thm:main}.]
    Let $C=20^4$ and $n>C^{10}$. Define $\alpha_n = 11C/\sqrt[4]{n}$, and suppose $\chi \in \X_n$ is such that $f(n,\chi)>\sqrt{n}$. By \Cref{prop:main}, we have $f(m)<\sqrt{m} +C$ for every $m \in \N$. Applying \Cref{lemma:long-path} with $C_1=C$ and $C_2=0$, we see that without loss of generality, there exists a blue path $P$ such that the size of $Y=[n]\setminus V(P)$ is at most $(1+\alpha_n)\sqrt{n}$ and with the additional property that every vertex in $Y$ has at most $2(C+1)\sqrt{n}$ blue neighbours in $X =V(P)$. By \Cref{lemma:Y-bound}, if $\vert Y_0\vert \leq (1-2\alpha_n)\sqrt{n}$ or $\vert Y\vert\leq \sqrt{n}-1$, then $f(n,\chi)\leq\sqrt{n}$, a contradiction. Let us therefore assume that $\vert Y_0\vert >(1-2\alpha_n)\sqrt{n}$ and $\vert Y\vert > \sqrt{n}-1$. 

    We define $H$ to be the bipartite graph on vertex set $X\cup Y$ whose edges are the red edges between $X$ and $Y$ under $\chi$. Let $Y_1 = Y \setminus Y_0$, $X_0 = \{x\in X: d_H(x)=|Y|\}$, and $X_1 = X\setminus X_0$. We have $|Y_1|=|Y|-|Y_0| <(1+\alpha_n)\sqrt{n}-(1-2\alpha_n)\sqrt{n}  = 33Cn^{1/4}$. By choice of $P$, we know that $\vert \{x\in X: \chi(xy) = \text{blue}\}\vert < 2(C+1)\sqrt{n}$ for every $y\in Y$, which implies that $\vert X_1 \vert < 2(C+1)\sqrt{n}|Y_1|< 67C^2n^{3/4}$. 
    Now, either $Y_1 = \emptyset = X_1$, or, since $n$ is sufficiently large, we have
    \begin{equation*}
        \frac{|X_0|}{|Y_1|}= \frac{|X|-|X_1|}{|Y_1|}>\frac{n-(1+\alpha_n)n^{1/2}-67C^2n^{3/4}}{33C n^{1/4}}>n^{1/2}>\frac{134C^2n^{3/4}}{(1-2\alpha_n) n^{1/2}}>\frac{2|X_1|}{|Y_0|}.
    \end{equation*}
Either way, we may apply \Cref{lemma:bipartite-decomposition-Y0} to obtain a collection of at most $\lceil |X|/(|Y|+1) \rceil$ paths in $H$ covering all vertices in $H$. This corresponds to a red path cover in $K_n$, which implies that $f(n,\chi)\leq \lceil |X|/(|Y|+1) \rceil$. Since $|Y|>\sqrt{n} -1$, it follows that
\begin{equation*}
    f(n,\chi)\leq \left\lceil \frac{|X|}{|Y|+1}\right \rceil = \left\lceil \frac{n+1}{|Y|+1}\right \rceil -1
    \leq \left\lceil \sqrt{n+1} \right\rceil -1
   \leq \sqrt{n}.
\end{equation*}
Thus we have a contradiction, proving the theorem.
\end{proof}

\section*{Acknowledgements}
The second author is supported by the London Mathematical Society and the Heilbronn Institute for Mathematical Research through an Early Career Fellowship. The third author is supported by the Martingale Foundation. We thank Zach Hunter and Frank Oberson for their comments on the paper.
\bibliographystyle{plain}
\bibliography{bib}

\end{document}